\documentclass[10pt]{article}

\usepackage{amssymb}
\usepackage{amsfonts}
\usepackage{amsmath}
\usepackage{amscd}
\usepackage[all,cmtip]{xy}
\usepackage{amsthm}
\usepackage{setspace}
\usepackage{enumerate}
\usepackage{geometry}
\usepackage{tcolorbox}
\usepackage{hyperref}

\hypersetup{
allbordercolors=blue, pdfborder=0 0 1.3}

\theoremstyle{plain}
\newtheorem{thm}{Theorem}

\newtheorem{cor}[subsection]{Corollary}

\theoremstyle{definition}

\title{Kirillov polynomials for the exceptional Lie algebra $\mathfrak g_{2}$}
\author{Martin T. Luu \footnote{Department of Mathematics, University of California, Davis, email: mluu@math.ucdavis.edu}} 

\date{}

\begin{document}

\maketitle

\begin{abstract}
As part of the development of the orbit method, Kirillov has counted the number of strictly upper triangular matrices with coefficients in a finite field of $q$ elements and fixed Jordan type. One obtains polynomials with respect to $q$ with many interesting properties and close relation to type A representation theory. In the present work we develop the corresponding theory for the exceptional Lie algebra $\mathfrak g_2$. In particular, we show that the leading coefficient can be expressed in terms of the Springer correspondence.
\end{abstract}

\section{Introduction}

Let $\mathbb{F}_{q}$ be a finite field of characteristic $p$ and consider the set $\mathfrak n(n,\mathbb{F}_{q})$ of strictly upper-triangular $n\times n$ matrices with entries in $\mathbb{F}_{q}$. The group $\mathcal G_{n}(\mathbb{F}_{q})$ of upper-triangular $n\times n$ matrices with $1$'s on the diagonal acts on $\mathfrak n(n,\mathbb{F}_{q})$ via conjugation. The (complicated) structure of the adjoint orbits plays an important role in the orbit method in representation theory. As an approximation to the structure of the adjoint orbits, Kirillov initiated in \cite{KIR1}, \cite{KIR2} a detailed study of the number $P_{\lambda}(q)$ of elements in $\mathfrak n(n,\mathbb{F}_{q})$ of fixed Jordan type $\lambda$ (here $\lambda$ is a partition of $n$). One obtains polynomials with respect to $q$ (we will refer to them as Kirillov polynomials) with many interesting properties that are strongly influenced by type A representation theory. For example, the leading coefficient is given by the dimension of the irreducible representation $V_{\lambda}$ of the permutation group $S_{n}$ associated to $\lambda$. 

To illustrate the huge simplification that occurs when passing from adjoint orbits to Jordan types let us consider the case $n=4$. Using a recursion relation with respect to $n$, one can calculate, see \cite{FK}, that for $n=4$ the Kirillov polynomials are
\begin{eqnarray*}
P_{4}(q)&=&q^3 \cdot (q-1)^3\\
P_{3,1}(q)&=&q^2  \cdot (q-1)^2  \cdot(1+3q)\\
P_{2,2}(q)&=&q \cdot (q-1)^2 \cdot (1+2q)\\
P_{2,1,1}(q)&=&(q-1) \cdot(1+2q+3q^2)\\
P_{1,1,1,1}(q)&=&1
\end{eqnarray*}

The elements in $\mathfrak n(4,\mathbb{F}_{q})$ correspond bijectively to those in $\mathcal G_4(\mathbb{F}_{q})$ via $X \mapsto \textrm{id}+X$ and this maps adjoint orbits to conjugacy classes. The $\mathcal G_4(\mathbb{F}_{q})$ conjugacy classes are classified in \cite{VLA1} in the following manner. Every conjugacy class contains a unique so-called primitive element and these are classified by their type: Each of the $6$ possibly non-zero entries is either a ramification point with $0$ entry (denoted by $\theta$), a ramification point with non-zero entry (denoted by $\bullet$), or an inert point with $0$ entry (denoted by $0$). We refer to \cite{VLA1} for the precise definitions. Order the indices containing potentially non-zero matrix entries as
$$(3,4)<(2,3)<(2,4)<(1,2)<(1,3)<(1,4)$$
With respect to this ordering one lists the types of each matrix entry, for example  matrices of type $\theta, \theta,\bullet, \theta, \theta ,0$ are of the shape
$$\begin{bmatrix}
1& \theta & \theta & 0\\
0&1 &\theta &\bullet \\
0&0&1&\theta\\
0&0&0&1
\end{bmatrix}$$
The classification of conjugacy classes in $\mathcal G_n(\mathbb{F}_{q})$ (for $n$ sufficiently small) in terms of the types of primitive elements is achieved in \cite{VLA1}, \cite{VLA2} and subsequent papers. To use these results (which are available only for small $n$) to recover the Kirillov polynomials it remains to calculate the Jordan type for each conjugacy class. For example, the conjugacy classes corresponding to the type $\theta, \theta,\bullet, \theta, \theta ,0$ yield $q\cdot (q-1)$ matrices of Jordan type corresponding to the partition $4=2+1+1$. This holds since the type has one $\bullet$ entry yielding $q-1$ choices and the size of the centralizer of a canonical matrix of this type is given in \cite {VLA1} as $q^5$ and therefore the size of the conjugacy class is $q^{6-5}$. We list the analogous calculations for all conjugacy classes:
 
 \hspace{0.2in}
 
\begin{center}

\begin{tabular}{c|c|c}
\textrm{conjugacy type} & \textrm{Jordan type} & \textrm{number} \\ [2pt]
\hline
\hline 
&&\\
$\theta,\theta,\theta,\theta,\theta,\theta$ & 1,1,1,1 & 1  \\ [2pt] 
$\theta,\theta,\theta,\theta,\theta,\bullet$ & 2,1,1& $q-1$ \\[2pt] 
$\theta,\theta,\bullet, \theta,\theta,0 $& 2,1,1& $q(q-1)$ \\[2pt]
$\theta,\bullet,0, \theta,0, \theta $& 2,1,1& $q^2(q-1)$\\[2pt]
$\bullet, \theta, 0, \theta,\theta,0 $& 2,1,1 & $q^2(q-1)$\\[2pt]
$\bullet, \theta, 0 ,\bullet, \bullet, 0 $& 3,1 & $q^2(q-1)^3 $\\[2pt]
$\theta,\theta,\theta,\theta,\bullet,0 $& 2,1,1 & $q(q-1)$\\[2pt]
$\theta,\theta,\bullet,\theta,\bullet,0$ & 2,2 & $q(q-1)^2$\\[2pt]
$\theta,\bullet, 0,\theta,0,\bullet $& 2,2 & $q^2(q-1)^2$\\[2pt]
$\bullet, \theta, 0 , \theta, \bullet, 0 $& 3,1 & $q^2(q-1)^2$ \\[2pt]
$\bullet,\bullet, 0, \theta, 0, 0 $& 3,1 & $q^3 (q-1)^2$ \\[2pt]
$\theta,\theta,\theta, \bullet, 0,0 $& 2,1,1 & $q^2(q-1) $\\[2pt]
$\theta,\theta, \bullet,\bullet, 0,0 $& 3,1 & $q^2 (q-1)^2$ \\[2pt]
$\theta, \bullet, 0, \bullet, 0,0 $& 3,1 & $q^3(q-1)^2$ \\[2pt]
$\bullet, \theta, 0,\bullet, \theta, 0 $& 2,2&$ q^2(q-1)^2$\\[2pt]
$\bullet,\bullet,,0,\bullet, 0,0 $& 4& $q^3(q-1)^3$ \\
&&
\end{tabular}

\end{center}
Adding things up one obtains, as expected
\begin{eqnarray*}
P_4(q)&=&q^3  \cdot (q-1)^3\\
P_{3,1}(q)&=&q^2  \cdot (q-1)^3+q^2  \cdot (q-1)^2+q^3  \cdot (q-1)^2+q^2  \cdot (q-1)^2+q^3  \cdot (q-1)^2\\
&=&q^2 \cdot(q-1)^2  \cdot (1+3q)\\
P_{2,2}(q)&=&q  \cdot (q-1)^2  \cdot (1+q+q)\\
&=&q  \cdot (q-1)^2  \cdot (1+2q)\\
P_{2,1,1}(q)&=&(q-1)  \cdot (1+q+q^2+q^2+q+q^2)\\
&=&(q-1)  \cdot (1+2q+3q^2)\\
P_{1,1,1,1}(q)&=&1
\end{eqnarray*}
Note that the number of adjoint orbits for $n=4$ is $2q^3+q^2-2q$ but even the question whether for fixed $n$ the number of orbits is always a polynomial with respect to $q$ is still an open question. As $n$ grows, the classification of adjoint orbits quickly becomes unknown, whereas the Kirillov polynomials have a simple recursion relation that we now recall.

We write partitions as $\lambda=(\lambda_1,\cdots,\lambda_{N})$ with $\lambda_{i}$ non-increasing. Given $\lambda$ we denote by $\lambda'=(\lambda_1',\cdots, \lambda_{N'})$ its dual partition. A cell of the Young diagram of $\lambda$ is called removable if after its removal the Young diagram remains the Young diagram of a partition. Let $s$ denote the number of removable cells and let $(x_{j},y_{j})$ be the coordinates of the $j$'th removable cell, see \cite{FK} for coordinate conventions. Denote by $\lambda \downarrow j$ the partition obtained by removing the $j$'th removable cell. Then by loc. cit. (Proposition 3.1), one has 
$$P_{\lambda}(q)= \sum_{j=1}^{s} (q^{n-\lambda'_{y_{j}} }-q^{n-1-\lambda'_{y_{j}-1}} ) \cdot P_{\lambda \downarrow j}(q)$$
where we use the convention $\lambda_{0}'=\infty$.

This recursion allows the efficient calculation of the polynomials $P_{\lambda}(q)$, see \cite{FK} for a table of all polynomials with $n \le 10$. Many fundamental properties of the Kirillov polynomials can be obtained easily from the recursion. In particular, the occurrence of factors of the form $q$ and $q-1$ can be explained in this manner. One obtains
$$P_{\lambda}(q)=q^{\binom{n}{2}-\binom{N}{2}-\sum_{i=1}^{N'-1} \lambda_{i}' \lambda_{i+1}'}\cdot (q-1)^{n-N} \cdot R_{\lambda}(q)$$
with $R_{\lambda}(0)$ and $R_{\lambda}(1)$ non-zero. Furthermore, the polynomial $R_{\lambda}(q)=:\sum_{i} r_{i} q^i$ satisfies interesting properties: The constant term $r_0$ equals to $1$ and all coefficients are strictly positive integers, the degree of $R_{\lambda}(q)$ is expressed in term of the partition $\lambda'$ dual to $\lambda$ as
$$\textrm{deg } R_{\lambda}(q)= \sum_{i=1}^{N'-1} \lambda_{i}' \lambda_{i+1}'-\sum_{i=2}^{N'} \binom{\lambda'_{i}+1}{2}$$
Furthermore, the leading coefficient of $R_{\lambda}(q)$ is  $\textrm{dim } V_{\lambda}$ where as before $V_{\lambda}$ is the representation of the permutation group $S_n$ associated to $\lambda$. 

As an example consider
$$P_{3,2,1,1}(q)= -q^5 - 2 q^6 - 3 q^7 - 3 q^8 + 4 q^9 + 25 q^{10} + 11 q^{11} - 23 q^{12} - 43 q^{13} + 35 q^{14}$$
This factorizes as
$$P_{3,2,1,1}(q)=q^5\cdot (q-1)^3\cdot (1+ 5q+ 15q^2+ 34q^3+ 58q^4+ 62q^5+ 35q^6)$$
A natural question is if in addition to factors of the form $q$ and $q-1$ the Kirillov polynomials are typically divisible by other polynomials. For example the degree $6$ factor $R_{3,2,1,1}(q)$ in the above example turns out to be an irreducible polynomial in $\mathbb{Q}[q]$ (equivalently in $\mathbb{Z}[q]$). A first guess is that $R_{\lambda}(q)$ is always irreducible in the polynomial ring $\mathbb{Z}[q]$ but this is not quite true as the following factorizations show:
\begin{eqnarray*}
R_{3, 2, 1}(q) &=& (2q+1) \cdot (8q^3 + 8q^2 + 3q + 1)\\[2pt]
R_{4, 3, 1}(q) &=&(5q^2 + 4q + 1) \cdot (14q^3 + 10q^2 + 3q + 1)\\[2pt]
R_{5, 3, 1}(q) &=& (2q + 1) \cdot (81q^4 + 57q^3 + 23q^2 + 6q + 1)\\[2pt]
R_{4, 4, 1}(q) &=& (2q + 1) \cdot (42q^4 + 39q^3 + 18q^2 + 5q + 1)\\[2pt]
R_{4, 3, 2}(q) &=& (2q + 1) \cdot(84 q^8 + 195q^7 + 219q^6 + 171q^5 + 100q^4 + 47q^3 + 18q^2 + 5q + 1)\\[2pt]
R_{4, 2, 2, 1}(q) &=& (3q^2 + 2q + 1) \cdot(72q^6 + 111q^5 + 73q^4 + 38q^3 + 15q^2 + 5q + 1)\\[2pt]
R_{3, 3, 2, 1}(q) &=& (2q + 1) \cdot (84q^8 + 195q^7 + 219q^6 + 171q^5 + 100q^4 + 47q^3 + 18q^2 + 5q + 1)\\[2pt]
R_{7, 3}(q) &=& (5q + 1) \cdot (15q^2 + 4q + 1)\\[2pt]
R_{4, 4, 2}(q) &=&(2q + 1) \cdot (126q^6 + 180q^5 + 126q^4 + 62q^3 + 24q^2 + 6q + 1)
\end{eqnarray*}

However, among the more than a million partitions of all $n\le 50$ no more reducible cases show up when irreducibility is checked with SAGE.  So far we have not been able to rigorously establish irreducibility for those Kirillov polynomials not listed above. Partially motivated by such irreducibility questions, we establish in Theorem \ref{main-theorem} the existence and explicit formulas for the Kirillov polynomials for the exceptional Lie algebra $\mathfrak g_2$. The analogue of the $R_{\lambda}(q)$-polynomials for $\mathfrak g_{2}$ turn out to be all irreducible. 

There is a second phenomenon we address. As already indicated, the leading coefficient of $P_{\lambda}(q)$ and $R_{\lambda}(q)$ is given by the dimension of $V_{\lambda}$, the representation of the permutation group $S_n$ associated to $\lambda$. The partitions $\lambda$ can also be used to describe the nilpotent orbits via their Jordan normal forms. The association
$$\textrm{nilpotent orbit of type $\lambda$} \leadsto V_{\lambda}$$
is a special case of the Springer correspondence \cite{SPR1}, \cite{SPR2} that relates nilpotent orbits with representations of Weyl groups. Hence the leading order coefficients of the Kirillov polynomials are expressed in terms of the Springer correspondence. We show that this phenomenon persists for the exceptional Lie algebra $\mathfrak g_2$ (a general relation between Kirillov polynomials and the Springer correspondence will be described elsewhere).

\section{The case of $\mathfrak g_2$}

The analogue of the number of conjugacy classes in $\mathcal G_{n}(\mathbb{F}_{q})$ for Chevalley groups of type $G_2$ has been computed in \cite{GR} and is given by
$$q^3+2q^2-q-1$$
See also the related algorithm in \cite{BH}.  So, one can say that a study of Kirillov polynomials in this situation is not necessary, since the adjoint orbit structure is understood.  Nonetheless, this simple situation can be used to gain insight into the properties of Kirillov polynomials beyond type A.  In principle, one could approach the calculation of $\mathfrak g_2$ Kirillov polynomials by using the results of \cite{BH}, \cite{GR} and calculating Jordan types, with respect to a suitable representation, for representatives of each orbit. In the present work we carry out a more direct approach.

To define the notion of a Lie algebra of type $\mathfrak g_2$ over a finite field $\mathbb{F}_{q}$ one starts with a Chevalley basis of a complex simple Lie algebra of type $\mathfrak g_2$. Then consider the $\mathbb{F}_{q}$ vector space spanned by the Chevalley basis and view it as a Lie algebra by reducing the complex structure constants modulo $p$. As is customary, we assume from now on that the characteristic $p$ of $\mathbb{F}_{q}$ satisfies $p>3$. The next step for obtaining Kirillov polynomials is to choose a suitable representation. The paper \cite{HRT} is a useful reference for a general discussion of matrix representations of the complex exceptional Lie algebras. In particular, via folding $\mathfrak g_2$ from $\mathfrak s \mathfrak 0_{8}$ one obtains a faithful $7$-dimensional representation $\xi$. Let $\alpha_1,\alpha_2$ be a set of simple roots and assume $\alpha_1$ is the short root. The $6$ positive roots are then $\alpha_1,\alpha_2,\alpha_1+\alpha_2,2\alpha_1+\alpha_2, 3\alpha_1+\alpha_2, 3\alpha_1+2\alpha_2$. For a root $\alpha$ let $e_{\alpha}$ be the corresponding basis element of a Chevalley basis. A formula for $\xi(e_{\alpha_1})$ and $\xi(e_{\alpha_2})$ is given in \cite{HRT} (Section 3.6). Let $e_{i,j}$ be the $7\times 7$ matrix with $0$'s everywhere except a $1$ in the $i,j$ entry. Then
\begin{eqnarray}
\label{first-Chevalley-basis-equations}
\xi(e_{\alpha_1}) &=& e_{1,2}+2e_{3,4}+e_{4,5}+e_{6,7}\\
\label{second-Chevalley-basis-equations}
\xi(e_{\alpha_2}) &=& e_{2,3}+e_{5,6}
\end{eqnarray}
Possibly up to signs, the other Chevalley basis elements for positive roots are given by
\begin{eqnarray*}
e_{\alpha_1+\alpha_2} &=& [e_{\alpha_1},e_{\alpha_2}]\\
e_{2\alpha_1+\alpha_2}&=& \frac{1}{2} \cdot [e_{\alpha_1+\alpha_2},e_{\alpha_1}]\\
e_{3\alpha_1+\alpha_2} &=& \frac{1}{3} \cdot [e_{2\alpha_1+\alpha_2},e_{\alpha_1}]\\
e_{3\alpha_1+2\alpha_2}&=& [e_{3\alpha_1+\alpha_2},e_{\alpha_2}]
\end{eqnarray*}
Let
$$X:=\xi(a \cdot  e_{\alpha_1}+b \cdot  e_{\alpha_1+\alpha_2} +c \cdot e_{2\alpha_1+\alpha_2}+d \cdot e_{3\alpha_1+\alpha_2}+e \cdot e_{3\alpha_1+2\alpha_2}+ f \cdot e_{\alpha_2})$$ 
By Equation (\ref{first-Chevalley-basis-equations}) and Equation (\ref{second-Chevalley-basis-equations}) it follows that
$$X=\begin{bmatrix}
0& a& b& 2c& d& e& 0 \\
0& 0& f& -2b& -c& 0& e\\
0& 0& 0& 2a& 0& -c& -d\\
0& 0& 0& 0& a& b& c \\
0& 0& 0& 0& 0& f& -b\\
0& 0& 0& 0& 0& 0& a\\
0& 0& 0& 0& 0& 0& 0
\end{bmatrix}
$$

This yields a faithful representation of $\mathfrak g_{2}$ acting on $\mathbb{F}_{q}^{7}$. We show from first principles that there are indeed polynomials $P_{\lambda}(q)$ in $\mathbb{Z}[q]$ counting those elements $X$ of fixed Jordan type $\lambda$. We call them the $\mathfrak g_2$ Kirillov polynomials.

\begin{thm} 
\label{main-theorem}
Consider a finite field $\mathbb{F}_{q}$ of characteristic $p>3$.
The Kirillov polynomials $P_{\lambda}(q)$ for the exceptional Lie algebra $\mathfrak g_{2}$ with respect to its $7$-dimensional faithful representation exist and are given by
$$\begin{aligned}
P_{7}(q)&=q^{4}\cdot(q-1)^2\\[2pt]
P_{3,3,1}(q)&=q^2\cdot (q-1)^2\cdot (1+2q)\\[2pt]
P_{3,2,2}(q) &=q^2 \cdot (q-1) \cdot (1+2q)\\[2pt]
P_{2,2,1,1,1}(q)&=(q-1)\cdot (1+q+q^2)\\[2pt]
P_{1,1,1,1,1,1,1}(q)&=1
\end{aligned}$$
Write $$P_{\lambda}(q)=q^{a}\cdot (q-1)^{b} \cdot R_{\lambda}(q)$$ for $a,b$ in $\mathbb{Z}^{\ge 0}$ and $R_{\lambda}(0)R_{\lambda}(1)$ non-zero. The constant term of $R_{\lambda}(q)$ equals $1$, all other coefficients are strictly positive integers, and $R_{\lambda}(q)$ is irreducible in $\mathbb{Z}[q]$.
\end{thm}

By Theorem \ref{main-theorem} there are five Jordan types with non-zero $\mathfrak g_2$ Kirillov polynomial. This matches the number of nilpotent orbits in $\mathfrak g_2$ and for each such orbit the Springer correspondence \cite{SPR1}, \cite{SPR2} associates a representation of the Weyl group. We show that this correspondence yields the leading coefficients of the Kirillov polynomials:  

\begin{cor}
\label{Springer-corollary}
For each partition $\lambda$ of $7$ such that the $\mathfrak g_2$ Kirillov polynomial $P_{\lambda}(q)$ is non-zero there is a unique nilpotent orbit in $\mathfrak g_2$ that has Jordan type $\lambda$ in the $7$-dimensional representation. Let $V_{\lambda}$ be the complex
representation of the Weyl group of $\mathfrak g_{2}$ associated to this nilpotent orbit via the Springer correspondence. Then
$$P_{\lambda}(q)= (\textrm{\emph{dim }} V_{\lambda}) \cdot  q^{\textrm{\emph{deg }} P_{\lambda}(q)} + \textrm{\emph{ lower order terms}}$$
\end{cor}

We prove the theorem by calculating, in terms of the parameters $a, b, c, d, e, f$, the rank sequence $(r_1,\cdots, r_6)$ where  $r_i = \textrm{rank } X^{i}$. One sees that

$$ X^2=\begin{bmatrix}
0&0&af&0&ac&bc+df&2ae-2bd+2{c}^{
2}\\0&0&0&2af&-2ab&-2{b}^{2}-2cf&-bc-df\\ 
0&0&0&0&2{a}^{2}&2ab&ac\\ 
0&0&0&0&0&af&0\\ 
0&0&0&0&0&0&af\\ 
0&0&0&0&0&0&0\\ 
0&0&0&0&0&0&0
\end{bmatrix}
\;\; , \;\; X^3= \begin{bmatrix}
0&0&0&2\,{a}^{2}f&0&0&0\\
 0&0&0&0&2\,{a}^{2}f&0&0
\\0&0&0&0&0&2\,{a}^{2}f&0\\ 0&0&0
&0&0&0&{a}^{2}f\\0&0&0&0&0&0&0\\ 0
&0&0&0&0&0&0\\ 
0&0&0&0&0&0&0
\end{bmatrix}
$$

$$ X^4= \begin{bmatrix}
0&0&0&0&2\,{a}^{3}f&2\,{a}^{2}bf&2\,{a
}^{2}cf\\0&0&0&0&0&2\,{a}^{2}{f}^{2}&-2\,{a}^{2}bf
\\0&0&0&0&0&0&2\,{a}^{3}f\\ 0&0&0
&0&0&0&0\\ 0&0&0&0&0&0&0\\ 0&0&0&0
&0&0&0\\ 0&0&0&0&0&0&0
\end{bmatrix}
\; \; , \;\; X^5= \begin{bmatrix}
0&0&0&0&0&2\,{a}^{3}{f}^{2}&0
\\ 0&0&0&0&0&0&2\,{a}^{3}{f}^{2}
\\ 0&0&0&0&0&0&0\\ 0&0&0&0&0&0&0
\\ 0&0&0&0&0&0&0\\ 0&0&0&0&0&0&0
\\ 0&0&0&0&0&0&0
\end{bmatrix}  $$

$$X^6 = \begin{bmatrix}
0&0&0&0&0&0&2\,{a}^{4}{f}^{2}
\\ 0&0&0&0&0&0&0\\ 0&0&0&0&0&0&0
\\ 0&0&0&0&0&0&0\\ 0&0&0&0&0&0&0
\\ 0&0&0&0&0&0&0
\end{bmatrix}
$$

Since the overall number of matrices is $q^6$ it suffices to obtain the matrix count for all but one type of rank sequence. We will skip the count for the rank sequence $(4,2,0,0,0,0)$. Furthermore, from the explicit formulas for the powers of $X$ it follows that unless $af\ne0$ the ranks of $X^{4},X^{5},X^{6}$ are all zero.

\subsection*{Case $1$: $af \ne 0$}

The direct calculations of all power $X^{i}$ shows that the rank sequence is always $(6,5,4,3,2,1)$. Letting $a,f$ range through non-zero numbers and $b,c,d,e$ arbitrary one obtains $(q-1)^2 \cdot q^4$ matrices.

\subsection*{Case $2$: $a=0$ and $f\ne 0$}

An explicit calculation shows that the rank of $X$ equals the rank of
$$\begin{bmatrix}
f&-2b &-c&0&e\\
0&2(b^2+cf)& bc+df&0& 0\\
0&0&0&f&-b\\
0&0&0&0&b^2+cf\\
0&0&0&0&-(bc+df)
\end{bmatrix}$$
Hence if 
\begin{eqnarray}
\label{case2-vanishing-equation}
b^2+cf = 0 =bc+df
\end{eqnarray}
then the rank of $X$ is $2$ and otherwise the rank is $4$. The rank of $X^2$ equals the rank of
$$\begin{bmatrix}
bc+df & 2(c^2-bd)\\
-2(b^2+cf) & -(bc+df)
\end{bmatrix}$$
Since
$$f(c^2-bd)=(b^2+cf)c-(bc+df)b$$
it follows that if Equation (\ref{case2-vanishing-equation}) holds, then the rank of $X^2$ is $0$. If not both of $b^2+cf$ and $bc+df$ vanish then the rank of $X^2$ equals to $1$ if
\begin{eqnarray}
\label{2112-Equation}
(bc+df)^2-4(b^2+cf)(c^2-bd)=0
\end{eqnarray}
and the rank equals $2$ otherwise.

To sum up, the possible rank sequences are $(4,2,0,0,0,0)$, $(4,1,0,0,0,0)$, and $(2,0,0,0,0,0)$. The latter occurs precisely when Equation (\ref{case2-vanishing-equation}) holds. Let us count solutions $a,b,c,d,e,f$ that also satisfy $a=0$ and $f \ne 0$. If $b=0$ then $c=0$ and $d=0$ and so one obtains $q \cdot (q-1)$ solutions. If $b \ne 0$ let $f\ne 0$ be arbitrary. This uniquely determines $c$ via the first equation and then $d$ via the second equation. So one gets $q\cdot (q-1)^2$ solutions. Together one gets
$$q^2 \cdot (q-1)$$
matrices with rank sequence $(2,0,0,0,0,0)$. As we have seen, the rank sequence $(4,1,0,0,0,0)$ occurs precisely when
\begin{enumerate}[(i)]
\item
not both of $b^2+cf=0$ and $bc+df=0$ hold
\item
$(bc+df)^2-4(b^2+cf)(c^2-bd)=0$
\end{enumerate}

We start by counting solutions to Equation (\ref{2112-Equation}), ignoring for now  the condition $f\ne0$. We consider two cases, depending on whether $d$ vanishes or not.

\subsubsection*{Case $d\ne 0$}

View Equation (\ref{2112-Equation}) as a quadratic equation for $f$. Looking at the discriminant one sees that given $(b,c,d)$ there are solutions for $f$ if and only if $(c^2-bd)^3$ is a square, or equivalently if and only if $c^2-bd$ is a square. If it is a square and equal to $0$ then there is precisely one $f$ solution, otherwise there are two distinct solutions for $f$.

Suppose now  $c^{2}-bd$ is a square, so there is $x$ with
\begin{eqnarray}
\label{b-d-equation}
bd=c^2-x^2
\end{eqnarray}
For each of the $(q+1)/2$ possible values of $x^2$ one can choose $c$ arbitrary and since $d\ne 0$ one obtains $q-1$ pairs $(b,d)$ solving Equation (\ref{b-d-equation}). Hence there are $q\cdot (q-1)^2/2$ triples $b,c,d$ such that $c^2-bd$ is a non-zero square and there are $q\cdot (q-1)$ triples $b,c,d$ such that $c^{2}-bd=0$. Since $a=0$ and $e$ is arbitrary, the number of choices of $a,b,c,d,e,f$ that solve Equation (\ref{2112-Equation}) is
$$q^2\cdot (q-1) + \frac{q^2\cdot(q-1)^2}{2} \cdot 2=q^3\cdot (q-1)$$

\subsubsection*{Case $d=0$}

Equation (\ref{2112-Equation}) now yields
\begin{eqnarray}
\label{2112-equation-d}
(bc)^2- 4(b^2+cf)c^2= 0
\end{eqnarray}
If $c=0$ this automatically holds and this gives $q^3$ solutions for $a,b,c,d,e,f$ coming from choosing $b,e,f$ arbitrarily. For fixed $c\ne 0$ we are solving$-3b^2=4cf$. Letting $b$ be arbitrary, this uniquely determines $f$. Since $e$ is arbitrary, one obtains $q^2\cdot(q-1) $ solutions of Equation (\ref{2112-equation-d}) with $c\ne 0$. 

Adding up the solutions with $d=0$ and $d\ne 0$ one obtains 
$$q^3  \cdot (q-1)+q^3+q^2  \cdot (q-1)$$
solutions to Equation (\ref{2112-Equation}). We have not yet imposed the condition $f \ne 0$. If $f=0$ we are solving
$$(bc)^2-4 b^2(c^2-bd)=0$$
If $b=0$ this holds automatically and one gets $q^3$ solutions for $a,b,c,d,e,f$. If $b\ne0$ we are solving $3c^2=4bd$ and in analogy with a previous calculation one obtains $q^2\cdot(q-1)$ solutions.

In conclusion, the amount of solutions to Equation (\ref{2112-Equation}) with $f \ne 0$ is
$$q^3  \cdot (q-1)$$
We want to count the number of these solutions that additionally satisfy $b^2+cf\ne 0$ or $bc+df\ne 0$.
As calculated earlier, the number of choices for $a,b,c,d,e,f$ with $a=0$, $f\ne 0$ and
$$b^2+cf =0 = bc+df$$
is given by $q^2  \cdot (q-1)$. It follows that the number of matrices with rank sequence $(4,1,0,0,0,0)$ is
$$q^3  \cdot (q-1)-q^2  \cdot (q-1)=q^2  \cdot (q-1)^2$$

\subsection*{Case $3$: $a\ne 0$ and $f=0$}

An explicit calculation shows that $X$ has the same rank as
$$\begin{bmatrix}
a & b & 2c & d & e & 0\\
0&0& 2a&0&-c&-d\\
0&0&0&a&b&c\\
0&0&0&0&0&a
\end{bmatrix}$$
and hence always has rank $4$. The rank of $X^2$ equals the rank of

$$\begin{bmatrix}
2a^2 & 2ab& ac\\
0 & 0 & 4ae-4bd +3c^2
\end{bmatrix}$$
and hence the rank of $X^2$ is $1$ if and only if
$$4ae-4bd+3c^2=0$$
and equals $2$ otherwise. Let us count the matrices with rank sequence $(4,1,0,0,0,0)$. Let $b,c,d$ be arbitrary. If $4bd=3c^2$ then $e=0$ and $a \ne 0$ is arbitrary, so one gets $q-1$ possible pairs $(a,e)$. If $4bd \ne 3c^2$ then $a\ne 0$ is arbitrary and uniquely determines $e$, so again one obtains $q-1$ possible pairs $(a,e)$. Hence, in total one obtains $q^3 \cdot (q-1)$ matrices. 

\subsection*{Case $4$: $a=0$ and $f=0$}

The rank of $X$ is $4$ if not both of $b,c$ are $0$. If both of $b,c$ are $0$ then the rank is $2$ if not both of $d,e$ are $0$. If $b=c=d=e=0$ then the rank is $0$. Let us now calculate the rank of $X^2$. If $b\ne 0$ then one sees that the rank is $1$ if $3c^2-4bd=0$ and the rank is $2$ otherwise. If $b=0$ then the rank is $0$ if $c=0$ and $1$ otherwise.

We now count how many matrices yield the rank sequence $(4,1,0,0,0,0)$. This happens exactly if one of the following two cases holds:

\begin{itemize}
\item $b\ne 0$ and $3c^2-4bd=0$
\item $b=0$ and $c \ne 0$
\end{itemize}
Consider the first case. If $c\ne 0$ then $d \ne 0$ is arbitrary and $b$ is uniquely determined, if $c=0$ then $d=0$ and $b\ne 0$ is arbitrary. The variable $e$ is always arbitrary. So one obtains
$$q \cdot((q-1)^2+(q-1))=q^2  \cdot  (q-1)$$
matrices. The second case yields $q^2  \cdot (q-1)$ matrices and hence in total one obtains
$$2 q^2  \cdot (q-1)$$
matrices with rank sequence $(4,1,0,0,0,0)$. The rank sequence $(2,0,0,0,0,0)$ occurs if $a=b=c=f=0$ and not both of $d,e$ are $0$. So this occurs for $(q-1) \cdot (q+1)$ matrices. Finally, the rank sequence $(0,0,0,0,0,0)$ occurs for one matrix.

\subsection{Completion of the proof}
Having completed the case by case analysis, we now count it all up. The rank sequence $(r_{1},\cdots, r_{6})$ yields the partition of $7$ in which $1\le i \le 7$ has multiplicity $\alpha_i$ given by
$$\alpha_i = r_{i-1}-2r_{i}+r_{i+1}$$
Here we put $r_{0}=\textrm{rank } X^{0}=7$ and $r_{i} =\textrm{rank }X^{i}=0$ for $i \ge 7$. So the rank sequences that showed up correspond to the following partitions:
\begin{eqnarray*}
(6,5,4,3,2,1)  &\leadsto & 7=7 \\
(4,2,0,0,0,0) &\leadsto & 7=3+3+1\\
(4,1,0,0,0,0) &\leadsto & 7= 3+2+2\\
(2,0,0,0,0,0) &\leadsto & 7= 2+2+1+1+1\\
(0,0,0,0,0,0) &\leadsto & 7= 1+1+1+1+1+1+1
\end{eqnarray*}
Adding up our various matrix counts one obtains
\begin{eqnarray*}
P_{7}(q) &=&  q^4\cdot (q-1)^2 \\[2pt]
P_{3,2,2}(q) &=& q^2\cdot (q-1)^2+q^3\cdot (q-1)+2q^2\cdot (q-1) = q^2\cdot (q-1) \cdot (1+2q)\\[2pt]
P_{2,2,1,1,1}(q) &=& q^2\cdot (q-1)+(q-1)\cdot (q+1)=(q-1)\cdot(1+q+q^2)\\[2pt]
P_{1,1,1,1,1,1,1}(q)&=&1
\end{eqnarray*}
As mentioned earlier, the count for one partition type comes for free if all others are known and we obtain
\begin{eqnarray*}
P_{3,3,1}(q)&=& q^6-q^4\cdot (q-1)^2  - q^2\cdot (q-1) \cdot (1+2q)-(q-1)\cdot(1+q+q^2) -1\\[2pt]
&=&q^2\cdot (q-1)^2\cdot (1+2q)
\end{eqnarray*}
This calculation of the Kirillov polynomials for $\mathfrak g_{2}$ completes the proof of Theorem \ref{main-theorem}. We now show that the leading coefficients of these polynomials can be expressed via the Springer correspondence.

The (complex) Springer correspondence \cite{SPR2} relates nilpotent orbits in the complex algebra $\mathfrak g_{2}$ with complex representations of the Weyl group of $\mathfrak g_{2}$. In fact, the nilpotent orbits themselves inject into the set of equivalence classes of representations and this injection is sufficient for our purposes. To obtain a bijection one enriches nilpotent orbits by representations of a suitable group (the original correspondence is obtained by choosing the trivial representations). For the explicit description of the Springer correspondence for $\mathfrak g_2$ see \cite{SPR1} (Section 7.16) and also \cite{CAR} (p. 427), \cite{SPR2}. What matters to us is the dimension of the representation associated to a nilpotent orbit and we record these in the table below. For each orbit we list a representative in terms of the previously described Chevalley basis and calculate the Jordan type, written as a partition of $7$, via the representation $\xi$. Since there are different terminological conventions for describing nilpotent orbits, we also list for each orbit the corresponding weighted Dynkin diagram: The weighted Dynkin diagram with $a$ corresponding to the short root and $b$ to the long root is denoted by $(a,b)$.

\begin{center}

\begin{tabular}{c|c|c|c|c}
\textrm{nilpotent orbit} & \textrm{weighted diagram} & \textrm{representative} & \textrm{partition} &\textrm{dimension } \\ 
\hline \hline
&&&&\\
$0$&& $0$ & $1+1+1+1+1+1+1$&$1$\\[2pt] 
$A_1$& $(0,1$)& $e_{3\alpha_1+2\alpha_2}$& $2+2+1+1+1$&$1$\\ [2pt] 
$\tilde A_1$ & $(1,0)$ & $e_{2\alpha_1+\alpha_2}$&$3+2+2$&$2$ \\ [2pt] 
$A_1+ \tilde A_1 $& $ (0,2)$ & $e_{\alpha_1}+e_{2\alpha_1+\alpha_2}$ &$3+3+1 $& $2$\\ [2pt] 
$G_2$ & $(2,2)$ & $e_{\alpha_1}+e_{\alpha_2} $& $7$& $1$\\
&&&&
\end{tabular}

\end{center}

By comparing the listed dimensions to the leading coefficients of the $\mathfrak g_2$ Kirillov polynomials one obtains Corollary \ref{Springer-corollary}.

\end{document}